\documentclass[12pt, oneside]{article}   	
\usepackage[margin=1.25in]{geometry}       
\usepackage{graphicx}				
\usepackage{amssymb}
\usepackage{amsmath}
\usepackage{tikz-cd}
\usepackage[mathcal]{eucal}
\usepackage{amsthm}
\theoremstyle{definition}
\newtheorem{thm}{Theorem}
\theoremstyle{definition}
\newtheorem{dfn}[thm]{Definition}
\theoremstyle{definition}
\newtheorem*{conv}{Conventions}
\theoremstyle{definition}

\theoremstyle{definition}

\theoremstyle{definition}
\newtheorem{lem}[thm]{Lemma}
\theoremstyle{definition}

\theoremstyle{definition}
\newtheorem{rem}[thm]{Remark}
\numberwithin{thm}{section}

\newcommand{\sub}{\subset}

\newcommand{\isom}{\approx}

\newcommand{\Rings}{\textup{Rings}}
\renewcommand{\O}{\mathcal{O}}

\newcommand{\m}{\mathfrak{m}}
\newcommand{\Sch}{\textup{Sch}}

\renewcommand{\lim}{\textup{lim}}
\renewcommand{\phi}{\varphi}
\newcommand{\id}{\textup{id}}

\usepackage[colorlinks,linkcolor=blue,citecolor=blue,urlcolor=blue,
            bookmarks,bookmarksnumbered]{hyperref}

\newcommand{\MTors}{\textup{$M$-Tors}}
\newcommand{\Exal}{\textup{Exal}}

\newcommand{\M}{\mathcal{M}}
\newcommand{\XSch}{X/\Sch}
\newcommand{\inv}{\textup{inv}}
\newcommand{\I}{\mathcal{I}}
\newcommand{\ThSch}{\textup{ThSch}}
\newcommand{\McoTors}{\textup{$\M$-coTors}}
\renewcommand{\O}{\mathcal{O}}

\title{First order thickenings and cotorsors}
\author{Nicholas Mertes\footnote{n.mertes@umiami.edu}\\ \textit{Department of Mathematics,}\\\textit{University of Miami},\\\textit{Coral Gables, FL}}

\date{June, 2020}							

\begin{document}

\maketitle

\begin{abstract}
Let $X$ be a scheme and let $\M$ be a quasi-coherent sheaf on $X$. Then $\M$ can be viewed as a cogroup object in the category of schemes under $X$. We show that the category of first order thickenings of $X$ by $\M$ is equivalent to the category of $\M$-cotorsors.
\end{abstract}

\section{Introduction}

\begin{conv}
We always use the word \textit{ring} to refer to a commutative ring with identity, and we always use the phrase \textit{ring homomorphism} to refer to an identity-preserving ring homomorphism. We write $\Rings$ for the category of rings and ring homomorphisms. We fix a ring $A$ and an $A$-module $M$. We write $\Rings/A$ for the category of rings over $A$.
\end{conv}

According to Beck~\cite{Beck}, we can view $M$ as a group object in $\Rings/A$ via the trivial square-zero extension $A\oplus M\to A$. Using this idea, the Andre-Quillen cohomology theory for rings~\cite{Quillen} was developed. The first Andre-Quillen cohomology group $H^1(A, M)$ can be identified with the set of isomorphism classes of square-zero extensions of $A$ by $M$. This cohomology theory is equivalent to a sheaf cohomology theory, and it is well-known that the first sheaf cohomology group can be interpreted as a set of isomorphism classes of sheaf-theoretic torsors. However, it would be more satisfying to interpret $H^1(A, M)$ in terms of torsors directly in the category $\Rings/A$ rather than in terms of sheaf-theoretic torsors. Given that $M$ can be viewed as a group object in $\Rings/A$, there should be a well-defined notion of $M$-torsors in $\Rings/A$. In section~\ref{sqandtors}, we define $M$-torsors and prove that the category of $M$-torsors is equivalent to the category of square-zero extensions of $A$ by $M$. This achieves our goal of interpreting $H^1(A, M)$ in terms of torsors directly in the category $\Rings/A$.

\begin{conv}
We write $\Sch$ for the category of schemes. We fix a scheme $X$ and a quasi-coherent sheaf $\M$ on $X$. We write $\XSch$ for the category of schemes under $X$.
\end{conv}

Given that $X$ can be viewed as a generalized ring and $\M$ can be viewed as a generalized module, it is reasonable to suspect that the above discussion can be generalized to the scheme-theoretic setting. Indeed, there is a scheme-theoretic analogue~\cite{Stacks} of Andre-Quillen cohomology which we will refer to as scheme-theoretic Andre-Quillen \textit{homology} because the morphisms in $\Sch$ are in the opposite direction of the morphisms in $\Rings$. The first scheme-theoretic Andre-Quillen homology group $H_1(X, \M)$ can be identified with the set of isomorphism classes of first order thickenings of $X$ by $\M$. One may wonder whether it is possible to interpret $H_1(X, \M)$ in terms of torsors. However, since $H_1(X,\M)$ is a homology group rather than a cohomology group, it is more reasonable to expect an interpretation in terms of \textit{cotorsors}. As shown in~\cite{Remy}, $\M$ can be viewed as a cogroup object in $\XSch$. Therefore, there should be a well-defined notion of $\M$-cotorsors in $\XSch$. In section~\ref{thickandcos}, we define $\M$-cotorsors and prove that the category of $\M$-cotorsors is equivalent to the category of first order thickenings of $X$ by $\M$. This achieves our goal of interpreting $H_1(X, \M)$ in terms of cotorsors.

\section{Square-zero extensions and torsors}\label{sqandtors}

The purpose of this section is to develop the commutative algebra which is necessary to give efficient descriptions of first order thickenings and cotorsors. A first order thickening is, roughly speaking, a morphism of schemes which gives rise to a square-zero extension in each stalk. Therefore, it is natural to start by discussing square-zero extensions. Throughout this section, we will only ever consider square-zero extensions of $A$ by $M$ and thus we will simply call these square-zero extensions.
\begin{dfn}
Let $f: B\to A$ be an object of $\Rings/A$. We write $f_\ast M$ for the $B$-module defined as follows. The underlying abelian group of $f_\ast M$ is the same as the underlying abelian group of $M$. The $B$-scalar multiplication on $f_\ast M$ is such that, for all $b\in B$ and $m\in f_\ast M$, $b\cdot m$ = $f(b)m$.
\end{dfn}
\begin{dfn}\label{defsquarezext}
A \textit{square-zero extension} is an object $f: B\to A$ of $\Rings/A$ such that $f$ is surjective and $(\ker(f))^2 = 0$, together with an isomorphism $\alpha_f: f_\ast M\to \ker(f)$ of $B$-modules.
\end{dfn}
Let $f: B\to A$ be a square-zero extension. Then we have an associated isomorphism $\alpha_f: f_\ast M\to \ker(f)$ of $B$-modules. Alternatively, we can define square-zero extensions in such a way that $\alpha_f$ is viewed as an isomorphism of $A$-modules rather than as an isomorphism of $B$-modules. In order to do this, we need to define an $A$-module structure on $\ker(f)$. Let $a\in A$ and let $b\in \ker(f)$. Since $f$ is surjective, choose $b'\in B$ such that $f(b') = a$. Then define the $A$-scalar multiplication on $\ker(f)$ such that $a\cdot b = b' b$. We need to check that this $A$-scalar multiplication is well-defined. Let $b''\in B$ be such that $f(b'') = a$. Then $b' - b''\in \ker(f)$. Since $(\ker(f))^2 = 0$, we have that
\[
\begin{split}
b'' b &= b''b + (b' - b'')b \\
&= b'b.
\end{split}
\]
Given this $A$-module structure on $\ker(f)$, it is straightforward to verify that $\alpha_f: M\to \ker(f)$ is an isomorphism of $A$-modules if and only if $\alpha_f: f_\ast M\to \ker(f)$ is an isomorphism of $B$-modules (recall that $M$ and $f_\ast M$ have the same underlying abelian group). In the literature on square-zero extensions, somehow $\alpha_f$ is unanimously viewed as an isomorphism of $A$-modules rather than as an isomorphism of $B$-modules. However, in this paper we always view $\alpha_f$ as an isomorphism of $B$-modules. This perspective becomes especially convenient in Section~\ref{thickandcos}.
\begin{dfn}\label{defmorphismsquarezero}
Let $f: B\to A$ and $g: C\to A$ be square-zero extensions. A \textit{morphism of square-zero extensions} from $f$ to $g$ is a morphism $h: f\to g$ in $\Rings/A$ such that, for each $m\in M$, $h(\alpha_f(m)) = \alpha_g(m)$. We write $\Exal(A, M)$ for the category of square-zero extensions and morphisms of square-zero extensions.
\end{dfn}
\begin{dfn}
We write $A\oplus M$ for the ring whose underlying abelian group is the direct sum of the abelian groups $A$ and $M$ and whose multiplication is such that, for all $(a_1, m_1), (a_2, m_2)\in A\oplus M$,
\[
(a_1, m_1)(a_2, m_2) = (a_1a_2, a_1m_2 + a_2m_1).
\]
The canonical projection $\pi_A: A\oplus M\to A$ is an object of $\Rings/A$ such that $\pi_A$ is surjective and $(\ker(\pi_A))^2 = 0$. We write $\alpha_{\pi_A}: (\pi_A)_\ast M\to \ker(\pi_A)$ for the isomorphism of $(A\oplus M)$-modules such that, for each $m\in M$, $\alpha_{\pi_A}(m) = (0, m)$. Then $\pi_A$ together with $\alpha_{\pi_A}$ is a square-zero extension.
\end{dfn}
We now want to equip $\pi_A: A\oplus M\to A$ with the structure of a group object in $\Rings/A$. In order to do this, we need a notion of products and a terminal object. Note that the identity $\id_A: A \to A$ is a terminal object in $\Rings/A$. Let $f: B\to A$ and $g: C\to A$ be objects of $\Rings/A$. We write $B\times_A C$ for the subring of $B\times C$ such that
\[
B\times_A C = \{(b, c)\in B\times C\,|\, f(b) = g(c)\}.
\]
We write $f\times g: B\times_A C\to A$ for the ring homomorphism such that, for each $(b, c)\in B\times_A C$, $(f\times g)(b, c) = f(b) = g(c)$.
\begin{dfn}
We equip $\pi_A: A\oplus M\to A$ with the structure of a group object in $\Rings/A$ as follows. We define $e_M: \id_A\to \pi_A$ such that, for each $a\in A$,
\[
e_M(a) = (a, 0).
\]
We define $+_M: \pi_A\times\pi_A\to\pi_A$ such that, for each $((a, m_1), (a, m_2))\in (A\oplus M)\times_A (A\oplus M)$,
\[
(a, m_1) +_M (a, m_2) = (a, m_1 + m_2).
\]
We define $\inv_M: \pi_A\to \pi_A$ such that, for each $(a, m)\in A\oplus M$,
\[
\inv_M(a, m) = (a, -m).
\]
\end{dfn}
It is straightforward to verify that $\pi_A$ together with the three morphisms $e_M$, $+_M$, and $\inv_M$ in $\Rings/A$ define a group object (in fact an abelian group object) in $\Rings/A$. Therefore, we can define $\pi_A$-torsors in $\Rings/A$. We will call these torsors $M$-torsors rather than $\pi_A$-torsors for convenience of notation.
\begin{dfn}
An \textit{$M$-torsor} is an object $f: B\to A$ of $\Rings/A$ such that $f$ is surjective, together with a morphism $\tau_f: \pi_A\times f\to f$ in $\Rings/A$ such that
\begin{enumerate}
\item \textbf{Convention}: For each $((a, m), b)\in (A\oplus M)\times_A B$, we will usually write $\tau_f(m, b)$ instead of $\tau_f((a, m), b)$. This should not cause confusion beacuse $((a, m), b)\in (A\oplus M)\times_A B$ implies that $a = f(b)$, and hence the $a$ is redundant in the notation.
\item For each $b\in B$, $\tau_f(0, b) = b$.
\item For each $(b_1, b_2)\in B\times_A B$, there exists a unique $m\in M$ such that
\[
\tau_f(m, b_2) = b_1.
\]
\end{enumerate}
\end{dfn}
One may suspect that we also need to require an associativity axiom. However, associativity is implied by our definition of $M$-torsor. We prove this associativity now, together with a preliminary lemma.
\begin{lem}
Let $f: B\to A$ be an $M$-torsor. If $m\in M$ and $b\in B$, then
\[
\tau_f(m, b) = \tau_f(m, 0) + b.
\]
\end{lem}
\begin{proof}
Let $m\in M$ and let $b\in B$. Then
\[
\begin{split}
\tau_f(m, b) &= \tau_f((f(b), m), b) \\
&= \tau_f((0, m) + (f(b), 0), b) \\
&= \tau_f(((0, m), 0) + ((f(b), 0), b)) \\
&= \tau_f((0, m), 0) + \tau_f((f(b), 0), b) \\
&= \tau_f(m, 0) + \tau_f(0, b) \\
&= \tau_f(m, 0) + b.
\end{split}
\]
\end{proof}
\begin{lem}
Let $f: B\to A$ be an $M$-torsor. If $m_1, m_2\in M$ and $b\in B$, then
\[
\tau_f(m_1 + m_2, b) = \tau_f(m_1, \tau_f(m_2, b)).
\]
\end{lem}
\begin{proof}
Let $m_1, m_2\in M$ and let $b\in B$. Then
\[
\begin{split}
\tau_f(m_1 + m_2, b) &= \tau_f(m_1 + m_2, 0) + b \\
&= \tau_f((0, m_1 + m_2), 0) + b \\
&= \tau_f((0, m_1) + (0, m_2), 0) + b \\
&= \tau_f(((0, m_1), 0) + ((0, m_2), 0)) + b \\
&= \tau_f((0, m_1), 0) + \tau_f((0, m_2), 0) + b \\
&= \tau_f(m_1, 0) + \tau_f(m_2, 0) + b \\
&= \tau_f(m_1, 0) + \tau_f(m_2, b) \\
&= \tau_f(m_1, \tau_f(m_2, b)).
\end{split}
\]
\end{proof}
\begin{dfn}
Let $f: B\to A$ and $g: C\to A$ be $M$-torsors. A \textit{morphism of $M$-torsors} from $f$ to $g$ is a morphism $h: f\to g$ in $\Rings/A$ such that, for all $m\in M$ and $b\in B$,
\[
h(\tau_f (m, b)) = \tau_g(m, h(b)).
\]
We write $\MTors$ for the category of $M$-torsors and morphisms of $M$-torsors.
\end{dfn}
\begin{dfn}\label{torstructonsqze}
Let $f: B\to A$ be a square-zero extension. We write $\tau_f: (A\oplus M)\times_A B\to B$ for the function such that, for all $m\in M$ and $b\in B$,
\[
\tau_f(m, b) = \alpha_f(m) + b.
\]
\end{dfn}
\begin{lem}\label{torstructonsqzelemma}
If $f: B\to A$ is a square-zero extension, then $\tau_f: (A\oplus M)\times_A B\to B$ gives $f$ the structure of an $M$-torsor.
\end{lem}
\begin{proof}
Let $f: B\to A$ be a square-zero extension. We first need to verify that $\tau_f: \pi_A\times f\to f$ is a morphism in $\Rings/A$. Let $((f(b_1), m_1), b_1), ((f(b_2), m_2), b_2)\in (A\oplus M)\times_A B$. Recall that $\alpha_f: f_\ast M\to \ker(f)$ is an isomorphism of $B$-modules, and recall that $(\ker(f))^2 = 0$. Then
\[
\begin{split}
\tau_f(((f(b_1), m_1), b_1)((f(b_2), m_2), b_2)) &= \tau_f((f(b_1), m_1)(f(b_2), m_2), b_1b_2)) \\
&= \tau_f((f(b_1)f(b_2), f(b_1)m_2 + f(b_2)m_1), b_1 b_2) \\
&= \tau_f(f(b_1)m_2 + f(b_2)m_1, b_1 b_2) \\
&= \alpha_f(f(b_1)m_2 + f(b_2)m_1) + b_1 b_2 \\
&= \alpha_f(b_1\cdot m_2 + b_2\cdot m_1) + b_1 b_2 \\
&= b_1 \alpha_f(m_2) + b_2 \alpha_f(m_1) + b_1 b_2 \\
&= \alpha_f(m_1) \alpha_f(m_2) + b_1 \alpha_f(m_2) + b_2 \alpha_f(m_1) + b_1 b_2 \\
&= (\alpha_f(m_1) + b_1)(\alpha_f(m_2) + b_2) \\
&= \tau_f(m_1, b_1)\tau_f(m_2, b_2) \\
&= \tau_f((f(b_1), m_1), b_1)\tau_f((f(b_2), m_2), b_2)
\end{split}
\]
and thus $\tau_f$ preserves multiplication. It is straightforward to show that $\tau_f$ preserves addition and preserves the identity, so $\tau_f$ is a ring homomorphism. Furthermore, $\tau_f$ is a morphism in $\Rings/A$ since, for all $m\in M$ and $b\in B$,
\[
\begin{split}
f(\tau_f(m, b)) &= f(\alpha_f(m) + b) \\
&= f(\alpha_f(m)) + f(b) \\
&= f(b) \\
&= (\pi_A\times f)(m, b). \\
\end{split}
\]
It remains to show that the morphism $\tau_f$ in $\Rings/A$ indeed satisfies the properties of an $M$-torsor. Note that, for each $b\in B$,
\[
\begin{split}
\tau_f(0, b) &= \alpha_f(0) + b\\
&= b.
\end{split}
\]
Now let $(b_1, b_2)\in B\times_A B$. Then $b_1 - b_2\in \ker(f)$. Since $\alpha_f: f_\ast M\to \ker(f)$ is bijective, there exists a unique $m\in M$ such that $\alpha_f(m) = b_1 - b_2$. Equivalently, there exists a unique $m\in M$ such that $\tau_f(m, b_2) = \alpha_f(m) + b_2 = b_1$.
\end{proof}
We now show that the assignment $\Psi: \Exal(A, M)\to \MTors$ which maps each square-zero extension to its corresponding $M$-torsor can be extended to a functor. In particular, we show that each morphism of square-zero extensions induces a corresponding morphism of $M$-torsors.
\begin{lem}\label{morphisquarezeroextimptormorph}
Let $f: B\to A$ and $g: C\to A$ be square-zero extensions. If $h: f\to g$ is a morphism of square-zero extensions, then $h$ is a morphism of $M$-torsors.
\end{lem}
\begin{proof}
Let $h: f\to g$ be a morphism of square-zero extensions. Let $m\in M$ and let $b\in B$. Then
\[
\begin{split}
h(\tau_f(m, b)) &= h(\alpha_f(m) + b) \\
&= h(\alpha_f(m)) + h(b) \\
&= \alpha_g(m) + h(b) \\
&= \tau_g(m, h(b)).
\end{split}
\]
\end{proof}
\begin{dfn}
We write $\Psi: \Exal(A, M)\to \MTors$ for the functor which maps each square-zero extension to its corresponding $M$-torsor and which maps each morphism of square-zero extensions to its corresponding morphism of $M$-torsors.
\end{dfn}
\begin{thm}\label{ringtorstheorem}
The functor $\Psi: \Exal(A, M)\to \MTors$ is an equivalence of categories.
\end{thm}
\begin{proof}
We will show that $\Psi$ is fully faithful and essentially surjective. Note that $\Psi$ is faithful since $\Psi$ maps each morphism of square-zero extensions to itself.

We now show that $\Psi$ is full. Let $f: B\to A$ and $g: C\to A$ be square-zero extensions and let $h: f\to g$ be a morphism of $M$-torsors. We want to show that $h$ is a morphism of square-zero extensions. Let $m\in M$. Then
\[
\begin{split}
h(\alpha_f(m)) &= h(\tau_f(m, 0)) \\
&= \tau_g(m, h(0)) \\
&= \tau_g(m, 0) \\
&= \alpha_g(m) + 0 \\
&= \alpha_g(m)
\end{split}
\]
and thus $h$ is a morphism of square-zero extensions.

Finally, we show that $\Psi$ is essentially surjective. Let $f: B\to A$ be an $M$-torsor. We want to show that $f$ can be given the structure of a square-zero extension in such a way that the $M$-torsor structure associated with the square-zero extension $f$ is the same as the original $M$-torsor structure on $f$.

First we show that $f: B\to A$ is such that $(\ker(f))^2 = 0$. Let $b_1, b_2\in \ker(f)$. Then $(0, b_1), (0, b_2)\in B\times_A B$, so there exist unique $m_1, m_2\in M$ such that $\tau_f(m_1, b_1) = 0$ and  $\tau_f(m_2, b_2) = 0$. Since $\tau_f: (A\oplus M)\times_A B\to B$ is a ring homomorphism, we have that
\[
\begin{split}
0 &= 0\cdot 0 \\
&= \tau_f(m_1, b_1)\tau_f(m_2, b_2) \\
&= \tau_f((0, m_1), b_1)\tau_f((0, m_2), b_2) \\
&= \tau_f(((0, m_1), b_1)((0, m_2), b_2)) \\
&= \tau_f((0, m_1)(0, m_2), b_1b_2) \\
&= \tau_f((0, 0), b_1b_2) \\
&= \tau_f(0, b_1b_2) \\
&= b_1b_2
\end{split}
\]
and thus $(\ker(f))^2 = 0$.

Note that, for each $m\in M$,
\[
\begin{split}
f(\tau_f(m, 0)) &= (f\circ\tau_f)((0, m), 0) \\
&= (\pi_A\times f)((0, m), 0) \\
&= 0
\end{split}
\]
and thus $\tau_f(m, 0)\in \ker(f)$. We define a function $\alpha_f: f_\ast M\to \ker(f)$ such that, for each $m\in M$, $\alpha_f(m) = \tau_f(m, 0)$. Then $\alpha_f$ is bijective since, for each $b\in \ker(f)$, $(b, 0)\in B\times_A B$ and thus there exists a unique $m\in M$ such that $\alpha_f(m) = \tau_f(m, 0) = b$. We now show that $\alpha_f: f_\ast M\to \ker(f)$ is an isomorphism of $B$-modules. Let $m_1, m_2\in M$. Then
\[
\begin{split}
\alpha_f(m_1 + m_2) &= \tau_f(m_1 + m_2, 0) \\
&= \tau_f(m_1, \tau_f(m_2,0)) \\
&= \tau_f(m_1, 0) + \tau_f(m_2, 0) \\
&= \alpha_f(m_1) + \alpha_f(m_2)
\end{split}
\]
and thus $\alpha_f$ is a group homomorphism. Let $b\in B$ and let $m\in M$. Then
\[
\begin{split}
\alpha_f(b\cdot m) &= \alpha_f(f(b)m) \\
&= \tau_f(f(b)m, 0) \\
&= \tau_f((0, f(b)m), 0) \\
&= \tau_f((f(b), 0)(0, m), 0) \\
&= \tau_f(((f(b), 0), b)((0, m), 0)) \\
&= \tau_f((f(b), 0), b)\tau_f((0, m), 0) \\
&= \tau_f(0, b)\tau_f(m, 0) \\
&= b\,\alpha_f(m)
\end{split}
\]
and thus $\alpha_f$ preserves $B$-scalar multiplication. Therefore, $\alpha_f: f_\ast M\to \ker(f)$ is an isomorphism of $B$-modules and hence $f$ together with $\alpha_f$ is a square-zero extension.

It remains to be shown that the $M$-torsor structure associated with the square-zero extension $f$ is the same as the original $M$-torsor structure on $f$. Let $m\in M$ and let $b\in B$. Then
\[
\begin{split}
\tau_f(m, b) &= \tau_f(m, 0) + b \\
&= \alpha_f(m) + b.
\end{split}
\]
Therefore, the functor $\Psi: \Exal(A, M)\to \MTors$ is an equivalence of categories.
\end{proof}

\section{First order thickenings and cotorsors}\label{thickandcos}

In effort to make this section easier to read, an attempt has been made to follow the organization of section~\ref{sqandtors} as closely as possible. In the previous section, we studied square-zero extensions. We will now describe the scheme-theoretic analogue of a square-zero extension, which is called a first order thickening. Since we will only ever consider first order thickenings of $X$ by $\M$, we will simply call these first order thickenings.
\begin{dfn}
Let $f: X\to Y$ be an object of $\XSch$ such that $f$ is a closed immersion. We write
\[
\I_f = \ker(f^\#: \O_Y\to f_\ast\O_X).
\]
\end{dfn}
Let $f: X\to Y$ be an object of $\XSch$ such that $f$ is a closed immersion and let $\alpha_f: f_\ast\M\to \I_f$ be an isomorphism of $\O_Y$-modules. Let $x\in X$. Then the canonical ring homomorphism $(f_\ast \O_X)_{f(x)}\to \O_{X, x}$ is an isomorphism and the canonical group homomorphism $(f_\ast \M)_{f(x)}\to \M_x$ is an isomorphism. Thus we obtain an object $f^\#_x: \O_{Y, f(x)}\to \O_{X, x}$ of $\Rings/\O_{X, x}$ such that $f^\#_x$ is surjective and an isomorphism $\alpha_{f, x}: (f^\#_x)_\ast \M_x\to \I_{f, f(x)}$ of $\O_{Y, f(x)}$-modules.
\begin{dfn}
See Definition~\ref{defsquarezext}. A \textit{first order thickening} is an object $f: X\to Y$ of $\XSch$ such that $f$ is a surjective closed immersion, together with an isomorphism $\alpha_f: f_\ast\M\to \I_f$ of $\O_Y$-modules such that, for each $x\in X$, the object $f^\#_x: \O_{Y, f(x)}\to \O_{X, x}$ of $\Rings/\O_{X, x}$ together with the isomorphism $\alpha_{f, x}: (f^\#_x)_\ast \M_x\to \I_{f, f(x)}$ of $\O_{Y, f(x)}$-modules is a square-zero extension of $\O_{X, x}$ by $\M_x$.
\end{dfn}
Let $f: X\to Y$ and $g: X\to Z$ be first order thickenings. Let $x\in X$. Then we obtain square-zero extensions $f^\#_x$ and $g^\#_x$ of $\O_{X, x}$ by $\M_x$. Let $h: f\to g$ be a morphism in $\XSch$. Then, for each $x\in X$, we obtain a morphism
\[
h^\#_x: \O_{Z, g(x)}\to \O_{Y, f(x)}
\]
in $\Rings/\O_{X, x}$.
\begin{dfn}
See Definition~\ref{defmorphismsquarezero}. Let $f: X\to Y$ and $g: X\to Z$ be first order thickenings. A \textit{morphism of first order thickenings} from $f$ to $g$ is a morphism $h: f\to g$ in $\XSch$ such that, for each $x\in X$,
\[
h^\#_x: \O_{Z, g(x)}\to \O_{Y, f(x)}
\]
is a morphism of square-zero extensions of $\O_{X, x}$ by $\M_x$. We write $\ThSch(X, \M)$ for the category of first order thickenings and morphisms of first order thickenings.
\end{dfn}
\begin{dfn}
We write $X\oplus \M$ for the scheme whose underlying topological space is $X$ and whose structure sheaf is $\O_{X\oplus \M} = \O_X\oplus \M$ where, for each open subset $U$ of $X$, $(\O_X\oplus\M)(U)$ is the ring $\O_X(U)\oplus \M(U)$. We write $i_X: X\to X\oplus \M$ for the morphism of schemes which is the identity on topological spaces and $i_X^\#: \O_X\oplus \M\to\O_X$ is the canonical projection of $\O_X\oplus \M$ onto $\O_X$. Then $i_X$ is a surjective closed immersion. Note that $(i_X)_\ast \M = \M$. We write $\alpha_{i_X}: \M\to \I_{i_X}$ for the isomorphism of $\O_{X\oplus M}$-modules such that, for all open subsets $U$ of $X$ and $m\in \M(U)$, $\alpha_{i_X}(U)(m) = (0, m)\in\O_X(U)\oplus \M(U)$. Then $i_X$ together with $\alpha_{i_X}$ is a first order thickening.
\end{dfn}
We now want to equip $i_X: X \to X\oplus \M$ with the structure of a cogroup object in $\XSch$. In order to do this, we need a notion of coproducts and an initial object. Note that the identity $\id_X: X \to X$ is an initial object in $\XSch$. Let $f: X\to Y$ and $g: X\to Z$ be objects of $\XSch$ such that $f$ and $g$ are closed immersions. First consider $f: X\to Y$ and $g:X\to Z$ as continuous maps of topological spaces. We can form the pushout
\[
\begin{tikzcd}
X \arrow{r}{g} \arrow[swap]{d}{f} & Z \arrow{d}{j_Z} \\
Y \arrow{r}{j_Y} & Y\amalg_X Z
\end{tikzcd}
\]
in the category of topological spaces. Then the topological space $Y\amalg_X Z$, together with the sheaf of rings
\[
\O_{Y\amalg_X Z} = (j_Y)_\ast\O_Y\times_{(j_Y\circ f)_\ast\O_X} (j_Z)_\ast\O_Z = (j_Y)_\ast\O_Y\times_{(j_Z\circ g)_\ast\O_X} (j_Z)_\ast\O_Z
\]
is a scheme and we write $f\amalg g: X\to Y\amalg_X Z$ for the associated object of $\XSch$. The object $f\amalg g$ of $\XSch$ is a coproduct of $f$ and $g$ since the scheme $Y\amalg_X Z$ is a pushout in the category of schemes~\cite{StacksPushout}. Thus we have a well-defined notion of cogroup object in $\XSch$, so long as the underlying ordinary object is a closed immersion.
\begin{dfn}
We equip $i_X: X\to X\oplus \M$ with the structure of a cogroup object in $\XSch$ as follows. We define $e_\M: i_X \to \id_X$ such that $e_\M$ is the identity on topological spaces and $e_\M^\#: \O_X\to \O_X\oplus \M$ is the canonical inclusion of $\O_X$ into $\O_X\oplus \M$. Note that, since $i_X$ is the identity on topological spaces, the coproduct $i_X\amalg i_X: X\to (X\oplus \M)\amalg_X (X\oplus\M)$ is also the identity on topological spaces. The structure sheaf of the scheme $(X\oplus \M)\amalg_X (X\oplus\M)$ is $(\O_X\oplus \M)\times_{\O_X} (\O_X\oplus\M)$. We define $+_\M: i_X\to i_X\amalg i_X$ such that $+_\M$ is the identity on topological spaces and $+_\M^\#: (\O_X\oplus \M)\times_{\O_X} (\O_X\oplus\M) \to \O_X\oplus \M$ is such that, for all open subsets $U$ of $X$ and $((a, m_1), (a, m_2))\in (\O_X(U)\oplus \M(U))\times_{\O_X(U)} (\O_X(U)\oplus\M(U))$,
\[
(a, m_1) +_{\M(U)}^\# (a, m_2) = (a, m_1 + m_2).
\]
We define $\inv_\M: i_X\to i_X$ such that $\inv_\M$ is the identity on topological spaces and $\inv_\M^\#: \O_X\oplus \M\to \O_X\oplus \M$ is such that, for all open subsets $U$ of $X$ and $(a, m)\in \O_X(U)\oplus \M(U)$,
\[
\inv_{\M(U)}^\#(a, m) = (a, -m).
\]
\end{dfn}

It is straightforward to verify that $i_X$ together with the three morphisms $e_\M$, $+_\M$, and $\inv_\M$ in $\XSch$ define a cogroup object (in fact an abelian cogroup object) in $\XSch$. Therefore, we can define $i_X$-cotorsors in $\XSch$. We will call these cotorsors $\M$-cotorsors rather than $i_X$-cotorsors for convenience of notation. Before stating the definition of $\M$-cotorsor, we need to make some preliminary remarks.

Let $f: X\to Y$ be an object of $\XSch$ such that $f$ is a closed immersion. Using the universal property of the coproduct $i_X\amalg f$, we obtain a morphism $\theta: (X\oplus\M)\amalg_X Y\to Y$ fitting into the pushout diagram
\[
\begin{tikzcd}
X \arrow[swap]{d}{i_X} \arrow{r}{f} &
Y \arrow{d}{j_Y} \arrow[ddr,bend left,"\id_Y"] \\
X\oplus \M \arrow{r}{j_{X\oplus\M}} \arrow[drr,bend right,"f\circ e_\M"'] &
(X\oplus\M)\amalg_X Y \arrow[dr,dashed,"\theta"] \\
&& Y
\end{tikzcd}
\]
Since $i_X$ is a cogroup object in $\XSch$, there should be a well-defined notion of a cogroup action of $i_X$ on $f$. A candidate for a cogroup action is a morphism $\tau_f: f\to i_X\amalg f$ in $\XSch$. For an ordinary group action on a space, we have an axiom which states that the identity element of the group acts trivially on all elements of the space. The analogue of that axiom in the cogroup setting states that $\theta\circ\tau_f = \id_Y$. Note that since the underlying topological space of the scheme $X\oplus \M$ is the topological space $X$, the topological space $Y$ is the underlying topological space of the scheme $(X\oplus \M)\amalg_X Y$. Therefore, at the level of topological spaces, the above pushout diagram becomes
\[
\begin{tikzcd}
X \arrow[swap]{d}{\id_X} \arrow{r}{f} &
Y \arrow{d}{\id_Y} \arrow[ddr,bend left,"\id_Y"] \\
X \arrow{r}{f} \arrow[drr,bend right,"f"'] &
Y \arrow[dr,dashed,"\theta"] \\
&& Y
\end{tikzcd}
\]
and we see that $\theta: Y\to Y$ must be the identity on topological spaces. Therefore, if $\tau_f: f\to i_X\amalg f$ is to be thought of as a cogroup action, we must have that $\theta\circ\tau_f = \id_Y$ and hence as functions on topological spaces we must have that $\tau_f = \theta = \id_Y$.

See Lemma~\ref{cotorimpsurj} and Remark~\ref{twodefscotors} in regards to the following definition.
\begin{dfn}
An \textit{$\M$-cotorsor} is an object $f: X\to Y$ of $\XSch$ such that $f$ is a closed immersion, together with a morphism $\tau_f: f\to i_X\amalg f$ in $\XSch$ such that
\begin{enumerate}
\item $\tau_f$ is the identity on topological spaces.
\item For each $y\in Y$, the object $f^\#_y: \O_{Y, y}\to (f_\ast \O_X)_y$ of $\Rings/(f_\ast \O_X)_y$ together with the morphism
\[
\tau^\#_{f,y}: ((f_\ast\O_X)_y\oplus (f_\ast \M)_y)\times_{(f_\ast \O_X)_y} \O_{Y, y}\to \O_{Y, y}
\]
in $\Rings/(f_\ast \O_X)_y$ is an $(f_\ast \M)_y$-torsor.
\end{enumerate}
\end{dfn}
A few comments need to be made about this definition of $\M$-cotorsor. Recall from the above discussion that the underlying topological space of $(X\oplus\M)\amalg_X Y$ is $Y$, and that $\tau_f$ must be the identity on topological spaces in order for $\tau_f$ to be thought of as a cogroup action. The structure sheaf of $(X\oplus\M)\amalg_X Y$ is $(f_\ast\O_X\oplus f_\ast\M)\times_{f_\ast \O_X} \O_Y$, and thus we have that $\tau_f^\#: (f_\ast\O_X\oplus f_\ast\M)\times_{f_\ast \O_X} \O_Y\to \O_Y$. Let $y\in Y$. Then we have a morphism
\[
\tau^\#_{f,y}: ((f_\ast\O_X)_y\oplus (f_\ast \M)_y)\times_{(f_\ast \O_X)_y} \O_{Y, y}\to \O_{Y, y}
\]
in $\Rings/(f_\ast \O_X)_y$. Since $f$ is a closed immersion, $f^\#_y: \O_{Y, y}\to (f_\ast \O_X)_y$ is surjective and thus it makes sense to test whether or not the object $f^\#_y$ of $\Rings/(f_\ast \O_X)_y$ together with $\tau^\#_{f, y}$ is an $(f_\ast \M)_y$-torsor.

Note that we do not require $f$ to be surjective in the definition of $\M$-cotorsor. However, the axioms for an $\M$-cotorsor imply that $f$ is surjective.
\begin{lem}\label{cotorimpsurj}
If $f: X\to Y$ is an $\M$-cotorsor, then $f$ is surjective.
\end{lem}
\begin{proof}
Let $f: X\to Y$ be an $\M$-cotorsor. Let $y\in Y$. Since $f$ is an $\M$-cotorsor, the object $f^\#_y: \O_{Y, y}\to (f_\ast \O_X)_y$ of $\Rings/(f_\ast \O_X)_y$ together with the morphism
\[
\tau^\#_{f,y}: ((f_\ast\O_X)_y\oplus (f_\ast \M)_y)\times_{(f_\ast \O_X)_y} \O_{Y, y}\to \O_{Y, y}
\]
in $\Rings/(f_\ast \O_X)_y$ is an $(f_\ast \M)_y$-torsor. Since every $(f_\ast \M)_y$-torsor has a square-zero kernel, we conclude that $(\ker(f^\#_y))^2 = 0$. Therefore, we must have that $\ker(f^\#_y)$ is a subset of the unique maximal ideal $\m_y\sub\O_{Y, y}$. Thus $\O_{Y, y}/\ker(f^\#_y)\isom (f_\ast \O_X)_y$ is nonzero and, since the image of $f$ is a closed subset of $Y$, we conclude that $y$ must be contained in the image of $f$.
\end{proof}
\begin{rem}\label{twodefscotors}
Given that the underlying object of $\XSch$ of any $\M$-cotorsor is automatically surjective, we obtain the following alternative definition of $\M$-cotorsor. Let $f: X\to Y$ be an object of $\XSch$ such that $f$ is a closed immersion and let $\tau_f: f\to i_X\amalg f$ be a morphism in $\XSch$. Then $f$ together with $\tau_f$ is an $\M$-cotorsor if and only if
\begin{enumerate}
\item $\tau_f$ is the identity on topological spaces.
\item $f$ is surjective.
\item For each $x\in X$, the object $f^\#_x: \O_{Y, f(x)}\to \O_{X, x}$ of $\Rings/\O_{X, x}$ together with the morphism
\[
\tau^\#_{f, x}: (\O_{X, x}\oplus \M_x)\times_{\O_{X, x}} \O_{Y, f(x)}\to \O_{Y, f(x)}
\]
in $\Rings/\O_{X, x}$ is an $\M_x$-torsor.
\end{enumerate}
This definition of $\M$-cotorsor is easier to work with than our original definition because its form is more closely related to first order thickenings. However, it seems to me that requiring $f$ to be surjective in the definition of $\M$-cotorsor is too strong of a condition.
\end{rem}
Let $f: X\to Y$ and $g: X\to Z$ be $\M$-cotorsors. Let $x\in X$. Then the objects $f^\#_x: \O_{Y, f(x)}\to \O_{X, x}$ and $g^\#_x: \O_{Z, g(x)}\to \O_{X, x}$ of $\Rings/\O_{X, x}$ together with the morphisms
\[
\tau^\#_{f, x}: (\O_{X, x}\oplus \M_x)\times_{\O_{X, x}} \O_{Y, f(x)}\to \O_{Y, f(x)}
\]
and
\[
\tau^\#_{g, x}: (\O_{X, x}\oplus \M_x)\times_{\O_{X, x}} \O_{Z, g(x)}\to \O_{Z, g(x)}
\]
in $\Rings/\O_{X, x}$ are $\M_x$-torsors.
\begin{dfn}
Let $f: X\to Y$ and $g: X\to Z$ be $\M$-cotorsors. A \textit{morphism of $\M$-cotorsors} from $f$ to $g$ is a morphism $h: f\to g$ in $\XSch$ such that, for each $x\in X$,
\[
h^\#_x: \O_{Z, g(x)}\to \O_{Y, f(x)}
\]
is a morphism of $\M_x$-torsors. We write $\McoTors$ for the category of $\M$-cotorsors and morphisms of $\M$-cotorsors.
\end{dfn}
\begin{dfn}
See Definition~\ref{torstructonsqze}, Lemma~\ref{torstructonsqzelemma}, and Remark~\ref{twodefscotors}. Let $f: X\to Y$ be a first order thickening. Then we define an $\M$-cotorsor structure on $f$ as follows. We write $\tau_f: f\to i_X\amalg f$ for the morphism in $\XSch$ such that $\tau_f$ is the identity on topological spaces and, for each $x\in X$,
\[
\tau^\#_{f, x}: (\O_{X, x}\oplus \M_x)\times_{\O_{X, x}} \O_{Y, f(x)}\to \O_{Y, f(x)}
\]
is the $\M_x$-torsor structure associated with the square-zero extension $f^\#_x: \O_{Y, f(x)}\to \O_{X, x}$ of $\O_{X, x}$ by $\M_x$.
\end{dfn}
In the above definition, it is straightforward to verify that the torsor structures specified on the stalks indeed give rise to a well-defined morphism $\tau_f: f\to i_X\amalg f$. We now show that the assignment $\Phi: \ThSch(X, \M)\to \McoTors$ which maps each first order thickening to its corresponding $\M$-cotorsor can be extended to a functor. In particular, we show that each morphism of first order thickenings induces a corresponding morphism of $\M$-cotorsors.
\begin{lem}
Let $f: X\to Y$ and $g: X\to Z$ be first order thickenings. If $h: f\to g$ is a morphism of first order thickenings, then $h$ is a morphism of $\M$-cotorsors.
\end{lem}
\begin{proof}
Let $h: f\to g$ be a morphism of first order thickenings and let $x\in X$. Then $h^\#_x: f^\#_x\to g^\#_x$ is a morphism of square-zero extensions of $\O_{X, x}$ by $\M_x$. Since $\tau^\#_{f, x}$ and $\tau^\#_{g, x}$ are the $\M_x$-torsor structures associated with the square-zero extensions $f^\#_x$ and $g^\#_x$ of $\O_{X, x}$ by $\M_x$, we know from Lemma~\ref{morphisquarezeroextimptormorph} that $h^\#_x$ is a morphism of $\M_x$-torsors.
\end{proof}
\begin{dfn}
We write $\Phi: \ThSch(X, \M)\to \McoTors$ for the functor which maps each first order thickening to its corresponding $\M$-cotorsor and which maps each morphism of first order thickenings to its corresponding morphism of $\M$-cotorsors.
\end{dfn}
\begin{thm}\label{mainthm}
The functor $\Phi: \ThSch(X, \M)\to \McoTors$ is an equivalence of categories.
\end{thm}
\begin{proof}
Now that we have set up all of the formalism for working with $\M$-cotorsors, this theorem follows straightforwardly from Theorem~\ref{ringtorstheorem}.
\end{proof}

\section*{Conclusion}

In this paper, we have shown that the first scheme-theoretic Andre-Quillen homology group (recall the discussion from the introduction) can be interpreted as a set of isomorphism classes of cotorsors. On a technical level, this result was proved by showing that first order thickenings of schemes can be identified with cotorsors for the cogroups associated with quasi-coherent sheaves. We have thus given a concrete geometric and categorical interpretation of an otherwise abstract homology group. Since the literature on cogroups and cotorsors in algebraic geometry is still in its infancy, this result marks an important turning point on the quest for a deeper understanding of homology theories in general.

We now make a few comments on how this research could be continued in the future. One project would be to take a fresh look at etale morphisms through the lens of cotorsors. The notion of a first order thickening is fundamentally related to etale morphisms of schemes, so it is reasonable to suspect that one may gain deeper insights into etale morphisms by replacing first order thickenings by their incarnation as cotorsors.

The next line of inquiry would be to study which other homology theories have convenient descriptions in terms of cotorsors, and also see which new homology theories can be developed from the cotorsor perspective. Finally, the cotorsor picture of the first scheme-theoretic Andre-Quillen homology group begs for a higher-categorical analogue. It is reasonable to suspect that a good theory of $\infty$-cotorsors should allow one to reproduce the full scheme-theoretic Andre-Quillen homology theory. Furthermore, there is nothing about this paper which is so special to schemes. There does not seem to be any barrier to proving similar results (perhaps with minor modifications) for algebraic spaces and stacks, as well as their derived analogues.

\section*{Acknowledgements}

In section 7.4.1 of Jacob Lurie's book \textit{Higher Algebra}~\cite{Lurie}, a comment is made about viewing square-zero extensions as torsors. That comment was the original source of motivation for this paper.

To the best of my knowledge, the only other reference to cogroups and cotorsors in the context of schemes is in the paper \textit{Combinatorial Differential Forms}~\cite{BreenMessing} by Lawrence Breen and William Messing. Some of what has been stated in the present paper can be found in section 1.11 of Breen and Messing's paper, but it does not appear that they carried out their inquiry all the way to the point of Theorem~\ref{mainthm}.

I would like to thank Sam Spiro for some helpful comments on the formatting of this paper.


\begin{thebibliography}\raggedright



\bibitem{Beck}
Triples, algebras and cohomology, Jonathan Mock Beck, Ph.D. thesis.\\
\url{http://www.tac.mta.ca/tac/reprints/articles/2/tr2abs.html}

\bibitem{Quillen}
Homology of commutative rings, Daniel Quillen, unpublished notes.\\
\url{https://en.wikipedia.org/wiki/Andre-Quillen_cohomology}

\bibitem{Stacks}
The Stacks project, 90: The Cotangent Complex.\\
\url{https://stacks.math.columbia.edu/tag/08P5}

\bibitem{Remy}
Automorphisms of categories of schemes, Remy van Dobben de Bruyn.\\
\url{https://arxiv.org/abs/1906.00921}

\bibitem{StacksPushout}
The Stacks project, 37.59: Pushouts in the category of schemes, II.\\
\url{https://stacks.math.columbia.edu/tag/0ECH}

\bibitem{Lurie}
Higher Algebra, Jacob Lurie.\\
\url{https://www.math.ias.edu/~lurie/}

\bibitem{BreenMessing}
Combinatorial Differential Forms, Lawrence Breen and William Messing.\\
\url{https://arxiv.org/abs/math/0005087}

\end{thebibliography}
\end{document}